\documentclass{article}
\usepackage[utf8x]{inputenc}
\usepackage[english]{babel}
\usepackage{graphicx}
\usepackage{amscd}
\usepackage{amsfonts}
\usepackage{latexsym}
\usepackage{amssymb,amsmath}
\usepackage[usenames]{color}

\usepackage{amscd}
\usepackage{psfrag}

\usepackage{amsthm}

\theoremstyle{plain}
\newtheorem{teor}{Theorem}[section]

\theoremstyle{definition}
\newtheorem{defi}[teor]{Definition}

\def\C{\mathbb{C}}
\def\S{\mathbb{S}}

\def\D{\mathbb{D}}

\def\l0{_{a_0}}
\def\i0{_{\iota_0}}

\begin{document}
\title{A note on parabolic-like maps}
\author{Luna Lomonaco}

\maketitle

\begin{abstract}
We show that the definition of parabolic-like map can be slightly modified, by asking $\partial \Delta$ to be a quasiarc out of the parabolic fixed point, instead of the dividing arcs to be $C^1$  on $[-1,0]$ and $[0,1]$. 
\end{abstract}
\section{Introduction}
Let $Per_1(1)$ denote the family of quadratic rational maps $P_A(z)=z+1/z+A,\,\,A\in \C$, with a parabolic fixed point of multiplier $1$ at $z= \infty$ and with critical points at $\pm 1$, modulo holomorphic conjugacy interchanging the critical points (this is, $P_A \sim_{conf} P_{-A}$).  
For $A \in \C$, we can define the filled Julia set $K_A$ of $P_A$ to be the complement of the basin of attraction of infinity (for $A=0$ the Julia set for map $P_0(z)=z+1/z$ is the imaginary axis, 
both the right and the left half planes are completely invariant parabolic basins of attractions, and we define $K_0$ to be the closure of the left half plane to be consistent with \cite{L1}).
A degree $2$ parabolic-like map is an object which behaves as a member $P_A$ of $Per_1(1)$ in a neighbourhood of its filled Julia set $K_A$. 
More generally, 
in \cite{L1} we defined parabolic-like maps in the following way:

\begin{defi}\label{definitionparlikemap} \textbf{(Parabolic-like maps)\,\,\,}
A \textit{parabolic-like map} of degree $d\geq2$ is a 4-tuple ($f,U',U,\gamma$) where 
\begin{itemize}
	\item $U'$ and $U$ are open subsets of $\C$, with $U',\,\, U$ and $U \cup U'$ isomorphic to a disc, and $U'$ not
contained in $U$,
	\item $f:U' \rightarrow U$ is a proper holomorphic map of degree $d\geq 2$ with a parabolic fixed point at $z=z_0$ of
 multiplier 1,
	\item $\gamma:[-1,1] \rightarrow \overline {U}$ is an arc with $\gamma(0)=z_0$, forward invariant under $f$, $C^1$
on $[-1,0]$ and on $[0,1]$, and such
that
$$f(\gamma(t))=\gamma(dt),\,\,\, \forall -\frac{1}{d} \leq t \leq \frac{1}{d},$$
$$\gamma([ \frac{1}{d}, 1)\cup (-1, -\frac{1}{d}]) \subseteq U \setminus U',\,\,\,\,\,\,\gamma(\pm 1) \in \partial U.$$
It resides in repelling petal(s) of $z_0$ and it divides $U'$ and $U$ into $\Omega', \Delta'$ and $\Omega, \Delta$
respectively, such that $\Omega' \subset \subset U$ 
(and $\Omega' \subset \Omega$), $f:\Delta' \rightarrow \Delta$ is an isomorphism %(see Fig. \ref{AAA}) 
and
$\Delta'$ contains at least one attracting fixed petal of $z_0$. We call the arc $\gamma$ a \textit{dividing arc}.

\end{itemize}

\end{defi}
In particular, we asked the dividing arc $\gamma:[-1,1] \rightarrow \overline {U}$ to be $C^1$
on $[-1,0]$ and on $[0,1]$. 

The goal of this note is to show that it is enough to ask the boundary of $\Delta$ to be a quasiarc outside the parabolic fixed point $z_0$, this is 
$$\partial \Delta \setminus \{z_0\} \mbox{is a quasiarc}.$$

\section{Motivation}\label{mot}
The definition \ref{definitionparlikemap} does not impose any restrictions for the boundaries $\partial U$ and $\partial U'$, and so it may seem counterintuitive to change this definition by relaxing the condition of the dividing arcs but imposing a condition on the whole $\partial \Delta$. 
The fact is that in the proof of  the Straightening Theorem of degree $2$ parabolic-like maps (see \cite{L1}) our first step is assuming $\partial U$ and $\partial U'$ to be $C^1$; and we can do it because a parabolic-like map is a local object, so if the boundaries are not $C^1$, we can replace them by subsets with $C^1$ boundaries (see \cite{L1}, first sentence of the proof of Theorem 6.3, page 26).
So, in definition \ref{definitionparlikemap}, \textit{de facto}, we want $\partial \Delta$ to be $C^1$ outside the parabolic fixed point $z_0$; and this is in order to prove that 
the boundary of the $\Delta$-part of the covering extension of the external map of the parabolic-like map
%the image of $\partial \Delta$ the Riemann map which under 
%the external map
is a \textit{quasicircle} (see \cite{L1}, second sentence of the proof of Claim 6.1, page 27, and for the definition of the covering extension of the external map of the parabolic-like map
see the end of Section \ref{covext}).  

%; and this suggests that we may directly replace the definition by asking a different condition which will ensure the image of $\partial \Delta$ under 
%the external map to be a \textit{quasicircle}, and this condition is
%$\partial \Delta \setminus \{z_0\}$  to be a quasiarc.\\

A quasiarc is the image of an interval under a quasiconformal map, while a $C^1$ curve is the image of an interval under a $C^1$ map. 
As diffeomorphisms on compact sets are quasiconformal (see \cite{Ah}), asking the dividing arc to be $C^1$ on  $[-1,0]$ and on $[0,1]$, and $\partial U$ to be $C^1$ implies that $\partial \Delta \setminus \{z_0\}$ is a quasiarc, and this together with Proposition 5.3(2) in \cite{L1} implies  that 
the boundary of the $\Delta$-part of the covering extension of the external map of the parabolic-like map is
%the image of $\partial \Delta$ under 
%the external map to be 
a \textit{quasicircle}. As this is exactly the condition needed in the proof of the Straightening Theorem, it suggests that we may directly replace the definition \ref{definitionparlikemap}
by asking $\partial \Delta \setminus \{z_0\}$  to be a quasiarc.\\

The practical reason behind this is that we may need to construct families of parabolic-like maps by holomorphic motion. 
In \cite{L2} we defined holomorphic families of parabolic-like maps to be families of parabolic-like maps defined as in \ref{definitionparlikemap}, and depending holomorphically on the parameter, and we proved that the connectedness locus of a nice proper family of parabolic-like map is homeomorphic to the parabolic Mandelbrot set without its root $M_1 \setminus \{0\}$, where a family is 'nice' if boundaries and the dividing arcs move holomorphically, and it is proper if it satisfies some compactness condition (see definition 5.3 in \cite{L2}). These conditions ensure the result about connectedness loci, but sometimes they are just too restrictive, in particular if one is in a setting for which it would be convenient to construct the family of parabolic-like maps by holomorphic motion.
A holomorphic motion is a basically a holomorphic isotopy. More precisely,
a holomorphic motion $\tau: \Lambda \times X \rightarrow \C$ with base point $\lambda_0 \in \Lambda$ of the set $X \subseteq \C$ is a map such that
\begin{enumerate}
\item for all $z \in X$, $\tau(\lambda_0,z)=z$;
\item if $z, w \in X,\,\, z \neq w$, then, for every $\lambda \in \Lambda$, $\tau(\lambda, z) \neq \tau (\lambda, w)$;
\item for every $z \in X$, $\tau(\lambda, z)$ is holomorphic in $\lambda$.
\end{enumerate}
By the $\lambda$-Lemma (\cite{mss}), if $\tau$ is a holomorphic motion, it is quasiconformal in $z$. Hence, if we move $\partial \Delta$ by holomorphic motion, we obtain that $\tau (\partial \Delta \setminus \{z_0\})$  is a quasiarc, even when $\partial \Delta \setminus \{z_0\}$ is $C^1$.

\section{About parabolic-like maps}
\subsection{Definitions}\label{covext}
All the  properties of parabolic-like maps are still valid after relaxing the definition. We summarize here the most important ones.
\begin{itemize}
\item The filled Julia set  $K_f$ of the parabolic-like map $(f,U',U,\gamma)$ is the set of points which never escape $\Omega'\cup{z_0}$ under iteration, and the Julia set $J_f$ is the boundary of the filled Julia set.
\item We say that $(f_1,U'_1,U_1, \gamma_1)$ is a parabolic-like restriction of $(f_2,U'_2,U_2, \gamma_2)$ if $U_1 \subseteq U_2$ and $f_1$ has the same degree and filled Julia set as $f_2$.
Two parabolic-like maps $(f,U_f',U_f,\gamma_f)$ and $(g,U_g',U_g,\gamma_g)$ are topologically (respectively quasiconformally or hybrid) equivalent if there exist parabolic-like restrictions
$(f,V_f',V_f,\gamma_f)$ and $(g,V_g',V_g,\gamma_g)$
homeomorphism $\phi:V_f \rightarrow V_g$ (respectively quasiconformal or quasiconformal with $\overline\partial \phi= 0$ a.e. on $K_f$) conjugating dynamics on $\Omega_{V_f}\cup \gamma_f$; while they are holomorphically conjugate if the homeomorphism $\phi$ between their parabolic-like restrictions is biholomorphic, $\phi(\gamma_f)=\gamma_g$, and it conjugates dinamics on the whoile $V_f'$. 
\item The external map of  the parabolic-like map $(f,U',U,\gamma)$ of degree $d$ is a degree $d$ real-analytic map $h_f: \S^1\rightarrow \S^1$ with a parabolic fixed point of multiplier $1$, and unique up to real-analytic differomorphism, which construction is done in the exact same way as in \cite{L1}. 
\end{itemize}
We recall the construction of the external map of a parabolic-like map  $(f,U',U,\gamma)$ with connected filled Julia set $K_f$, in order to show that the regularity of the dividing arc plays no role in this (and this is still the case if $K_f$ is disconnected, as the reader can check copying step by step the proof in \cite{L1}, pages 13-15).
Let $\alpha: \widehat \C \setminus K_f  \rightarrow  \widehat \C \setminus \overline \D$ be the Riemann map normalized fixing infinity with $\alpha(\gamma(t))\rightarrow 1$ as $t\rightarrow 0$,
let $W=\alpha(U)$, $W'=\alpha(U')$, and define
$$h_+:= \alpha\circ f \circ \alpha^{-1}: W' \rightarrow W.$$
Call $\widetilde W$ and $\widetilde W'$ the sets $W$ and $W'$ respectively union their reflections with respect to the $\S^1$ and union $\S^1$, reflect  $h_+$ with respect to $\S^1$ by the strong reflection principle, extend it to $\S^1$ and call $h: \S^1 \rightarrow \S^1$ the restriction to the unit circle, then by construction $h$ is a degree $d$ covering of the unit circle with a parabolic fixed point of multiplier $1$ at $z=1$, and we call it an external map for $f$. The external class $[h_f]$ of $f$ is the class of the map $h_f$ under conjugacy by real-analytic diffeomorphism.
Let $\gamma_h$ be the image under $\alpha$ of the dividing arc $\gamma$ for $(f,U',U,\gamma)$ (to be precise, $\alpha(\gamma)$ has two connected components, $\gamma^+_h$ and $\gamma^-_h$,
and $\gamma_h= \gamma^+_h \cup 1 \cup \gamma^-_h$), then
the covering extension of the external map of the parabolic-like map $(h,\widetilde W',\widetilde W,\gamma_h)$.

\subsection{The straightening theorem}
In \cite{L1} we proved the \textit{Straightening theorem for degree $2$ parabolic-like maps}: every degree $2$ parabolic-like map is hybrid equivalent to a member of $Per_1(1)$, a unique such a member if the filled Julia set is connected.  We obtained this theorem combining the following propositions and theorem in \cite{L1}:
\begin{itemize}
\item  For every $A \in \C$, the external class of $P_A$ is given by the class of $h_2(z)=\frac{z^2+1/3}{1+z^2/3}$ (Proposition 4.2)
\item A  degree $2$ parabolic-like map is holomorphically conjugate to a member of the family $Per_1(1)$ if and only if its external class is given by $h_2$ (Proposition 6.2)
\item Given a degree $d$ parabolic-like map $(f,U',U,\gamma_f)$ and a degree $d$ parabolic external map $h$, there exists a parabolic-like map $(g,V',V,\gamma_g)$ hybrid equivalent to $(f,U',U,\gamma_f)$ and with external map given by $[h]$ (Theorem 6.3).
\end{itemize}
The reader can check that dividing arcs play no role in the proofs of Proposition 4.2 and of Proposition 6.2, so we need to check that the proof of Theorem 6.3 in \cite{L1} is still valid.
\subsubsection{Proof of Theorem 6.3 in \cite{L1} is still valid}
As we said in Section \ref{mot}, the proof of Theorem 6.3 is based on proving that  the image of $\partial \Delta$ under 
the external map is a \textit{quasicircle} (see Claim 6.1), and this is needed in order to construct a quasiconformal map $\phi$ between the fundamental annulus $U_f \setminus \Omega'_f$ of the parabolic-like map $(f,U',U,\gamma)$ and a fundamental annulus of a parabolic-like restriction of a chosen external map $h$
%the map $h_2(z)=\frac{z^2+1/3}{1+z^2/3}$, which represents the external class of the family $P_A$. 
This quasiconformal map is then used in order to replace the standard complex structure on  the fundamental annulus $U_f \setminus \Omega'_f$ with the almost complex structure $\phi^*(\mu_0)$ which is the pullback under the quasiconformal map $\phi$ of the standard complex structure 'of' the map $h$. 
Then we spread by $f$ this new complex structure on $\Omega \setminus K_f$, and we basically leave it unchanged in $\Delta$, to obtain an almost complex structure on $U\setminus K_f$ invariant under $f$ (see Claim 6.2 in \cite{L1}), the integrating map of which conjugates by construction the parabolic-like map $(f,U',U,\gamma)$ with a parabolic-like map ($\phi \circ f \circ \phi^{-1}, \phi(U'),\phi(U), \phi(\gamma)$) of the same degree and with external class given by $h$.
So, if we started with a degree $2$ parabolic-like map, and we choose the external map to be $h_2$, the resulted parabolic-like map  
 ($\phi \circ f \circ \phi^{-1}, \phi(U'),\phi(U), \phi(\gamma)$) is  holomorphically conjugate  to a member of the family $Per_1(1)$.

The reader can check that in the proof of Claim 6.2 and in the rest of the proof of Theorem 6.3 in \cite{L1} the regularity of the dividing arcs play no role: their smoothness is  needed in the proof of Claim 6.1 to ensure that 
$\partial \Delta_{h_f}$ is a \textit{quasicircle}, where $h_f$ is the covering extension of the external map $h_f$ of $(f,U',U,\gamma_f)$. Replacing the $C^1$ condition on the dividing arcs with asking $\partial \Delta \setminus \{z_0\}$   to be a quasiarc, we obtain that $\partial \Delta_{h_f}\setminus \{1\}$ is a quasiarc. This quasiarc can be 'completed' in a quasicircle in the same way as in \cite{L1}: using Proposition 5.3(2), which show that, if dividing arcs of parabolic external maps $\gamma_{h_f,\pm}$, $\gamma_{h_2,\pm}$ are obtained as preimage of the same invariant curves under the respective Fatou coordinates $\phi_{h_f,\pm},\,\phi{h_2,\pm}$, then the composition of Fatou coordinates $\phi_{h_f}^{-1}\circ \phi_{h_f}: \gamma_{h_f} \rightarrow \gamma_{h_2}$ defined as
$$\phi_{h_f}^{-1}\circ \phi_{h_f}(z):= \left\{
\begin{array}{cl}
\phi_{h_f,+}^{-1}\circ \phi_{h_f,+}  &\mbox{ if } z \in \gamma_{h_f,+} \\
\phi_{h_f,-}^{-1}\circ \phi_{h_f,-}(z) &\mbox{ if  } z\in \gamma_{h_f,-} \\
\end{array}\right. 
$$
is quasisymmetric (see Proposition 5.3(2) in \cite{L1}).

\end{document}